\documentclass[preprint,12pt]{elsarticle}

\usepackage{graphicx}
\usepackage{amssymb}
\usepackage{amsmath}
\usepackage{lineno}
\usepackage{color}
\usepackage{tikz-cd}

\newsavebox{\smlmat}
\savebox{\smlmat}{$\left(\begin{smallmatrix}x&0\\0&y\end{smallmatrix}\right)$}

\journal{MAT437: K-theory and $C^*$-algebras}

\begin{document}

\begin{frontmatter}

\title{An Introduction to Abstract Classification Theory in the Operator Algebraic Setting}

\author{Patrick Fraser\\p.fraser@mail.utoronto.ca}

\address{University of Toronto}

\begin{abstract}

In the setting of modern mathematical logic and model theory, classification theory has been one of the landmark achievements of the field. Likewise, the classification of UHF-algebras and AF-algebras were substantial contributions to the field of operator algebra theory. These seemingly disparate topics of study in mathematics, model theory and operator algebras, have in recent years become closely related in many respects. I here attempt to bridge the gap between these two topics by discussing how operator algebraic classifications may be understood in terms of model-theoretic classification theory. This introductory article assumes basic familiarity with model theory and linear operator, but higher-level concepts are introduced when necessary. The focus of this introduction is conceptual and informal, and as such, many results are stated without proof, but relevant sources are cited for completeness. The reader should take this not as a detailed review, but rather as an overview of a general narrative thread connecting these two branches of modern mathematics.
\end{abstract}

\end{frontmatter}

\section{Introduction}

Freeman Dyson famously wrote about a distinction common to the sciences between the diversifiers and the unifiers~\cite{dyson_1985}. In order for a science to mature, it requires these two groups of people; the diversifiers who explore the vast jungle of the world, uncovering new, beautiful objects to study, proliferating novel ideas, and the unifiers, who fly above this jungle and piece together the distinct components, drawing them together into a coherent picture. Both are essential to understand the world fully. In mathematics, too, we see a vast array of diversifiers and unifiers. This essay focuses on the general goals and methods of the latter.

One of the common ways of unifying ideas in mathematics is through so-called classification. Initially, the mathematician may notice that there is a particularly interesting class of objects which have certain properties that may be used to understand how things work under a new light. Such was the case with the early formal introduction of what are now standard algebraic structures, such as groups, rings, and fields. Once such basic objects have been noticed, one then seeks to explore their behaviour. One might notice, for instance, that not all rings are commutative, but when they are, certain features appear which may be exploited to further our understanding. The initial piece of discovery opens a path, and mathematicians eagerly map their way along this new path, until they have a fairly detailed understanding of its cartography. Once the objects under consideration become sufficiently well understood, it becomes reasonable to try to ask the strongest, most general questions. Specifically, once mathematicians realize that a certain class of mathematical objects has interesting behaviour, they might want an explicit list of all of these objects of a particular form (up to some sort of equivalence). Often, obtaining such a complete list is very difficult; in many cases, a classification program first identifies a collection of quantities which are invariant under isomorphisms, and then claim that the objects in question are classified if these invariants are enough to tell apart any two non-isomorphic objects. The invariants themselves may be quite difficult to compute, and so a complete list is often unattainable (such is the case with von Neumann algebras, whose classification program is considered complete, even though it is often quite difficult to say whether or not two von Neumann algebras are isomorphic) but in principle, these invariants would allow one to compare any two objects and determine if they are isomorphic or not. This is the game of classification.

An ideal classification program would consist of a complete list, up to isomorphism, of all of the objects which behave in a certain way. In an elementary group theory course, one might, for instance see the historically significant classification of finite abelian groups. This classification tells us that a group $G$ is a finite abelian group if and only if $G\cong\mathbb{Z}/n_1\mathbb{Z}\times...\times \mathbb{Z}/n_k\mathbb{Z}$ for sufficient $n_i$'s and $k$; it gives us a complete list, up to isomorphism, of all possible finite abelian groups. In more sophisticated settings, classification theorems provide a powerful correspondence between the abstract and the concrete, and are usually landmark theorems in their disciplines.

Classification theory is a field of study in its own right, generally subsumed as a topic within model theory and mathematical logic. Here, I present an elementary introduction to classification theory in its model-theoretic setting, and then discuss how it may be concretely applied to the theory of operator algebras, UHF and AF-algebras in particular.

\section{Classification Theory}

I begin with the broad notions of classification theory from the model-theoretic standpoint, based on~\cite{shelah_2009}. In lieu of a basic introduction to model theory, I point the interested reader towards the very accessible introduction~\cite{leary_2015} and the slightly more sophisticated~\cite{marker_2002}. Assuming familiarity with the concept of a first-order predicate language,\footnote{Essentially, such a language is a meaningless assemblage of symbols which represent variables, constants, relations, and functions, with standard logical connectives like $\land$, $\rightarrow$, and $\lor$, as well as a negation symbol $\neg$ and universal and existential quantification symbols $\forall$ and $\exists$ respectively, which have basic rules for constructing grammatical terms, formulas, and sentences. The language itself is meaningless, but concrete meaning is prescribed to symbols and formulas when they are interpreted in a structure of that language.} I take the following definitions:\\

\textbf{Definition:} \textit{Fix a language $\mathcal{L}$. An $\mathcal{L}$-structure $\mathfrak{A}$ consists of a set of elements (the universe), a set of constants of the same cardinality as the constants of $\mathcal{L}$, and for every relation and function symbol in $\mathcal{L}$, a relation and function between elements of the universe.}\\

Languages have no meaning, but one may map the symbols of a language to a structure in which those symbols take on a meaning using an assignment function. Generally, one is concerned with \textit{models} of a particular set of formulas; that is, some structure in which all of a certain class of abstract formulas are true.\\

\textbf{Definition:} \textit{Let $\mathcal{L}$ be a language. Given a set $T$ of $\mathcal{L}$-formulas, an $\mathcal{L}$-structure $\mathcal{M}$ is a model of $T$ if $\mathcal{M}\vDash\phi$ for all $\phi\in T$.}\\

Structures are objects such as the natural numbers equipped with a binary addition function and a binary ordering relation, which may be taken as models for certain abstract formulas in the language. For example, in a language with one constant and an order relation $\mathcal{L}=\{0,\leq\}$, the abstract formula $\phi:=(\exists x)[(\forall y)(x\leq y)]$ is merely a string of symbols, but when interpreted in the structure of the natural numbers, this formula is satisfied because it means that the natural numbers have a smallest element, which is true. However, this same formula is not satisfied by the $\mathcal{L}$-structure that is the rational numbers (equipped with the usual ordering), for instance, because they do not have a lowest element under their usual ordering. As such, $\mathbb{N}$ models $\phi$, but $\mathbb{Q}$ does not (and we write $\mathbb{N}\vDash \phi$ and $\mathbb{Q}\not\vDash\phi$).

When we seek to ``classify" a collection of objects, we require first that they are in fact objects of the same kind; if we are concerned with groups, we are concerned with groups. If we are concerned with approximately finite $C^*$-algebras, we are concerned with approximately finite $C^*$-algebras. But we are never comparing apples with oranges. Classification takes place always within the confines of kind. Formally, this amounts to saying that we are restricting ourselves to $\mathcal{L}$-structures which are all grounded in the same language $\mathcal{L}$, and further, which all satisfy the same axioms. That is, they must be models of the same $\mathcal{L}$-theory. With this in mind, I give the following definitions which will help us, in broad strokes, to study the general behaviour of classes of models for a given theory:\\

\textbf{Definition:} \textit{Let $\mathfrak{K}=(K,\leq_\mathfrak{K})$ where $K=Mod(T)$ the set of all models of some first-order $\mathcal{L}$-theory $T$, and $\leq_\mathfrak{K}$ is partial ordering on $K$. Then $\mathfrak{K}$ is an abstract elementary class if: (i) $K$ is closed under isomorphism and $\leq_\mathfrak{K}$ preserves isomorphisms, (ii) if $\mathcal{M}\leq_\mathfrak{K}\mathcal{N}$, then $\mathcal{M}$ is a substructure of $\mathcal{N}$, (iii) (downward L\"owenheim-Skolem theorem) there exists a cardinal $\alpha\geq\aleph_0+|\mathcal{L}|$ such that for every $\mathcal{M}\in K$ and $A\subset |M|$, there is a $\mathcal{N}\in K$ with $\mathcal{N}\leq_\mathfrak{K}\mathcal{M}$, $A\subset |\mathcal{N}|$, and $||\mathcal{N}||\leq |A|+\alpha$, and (iv) (Tarski-Vaught Chains) $K$ is closed under $\leq_\mathfrak{K}$-increasing chains of arbitrary regular cardinality.}~\cite{shelah_2009, grossberg_2002}\\

Essentially, abstract elementary classes are classes of models for a particular theory which have enough particular constraints to ensure that the elements in each class all `look alike' in a precise sense. The models of a theory, particularly the infinite models, are in general not easy to pin down. Indeed, the following famous theorem ensures that, if there are any infinite models of some some first-order theory, there are always many distinct infinite models.\\

\textbf{Theorem (L\"owenheim-Skolem):} Suppose $T$ is a countable theory in a first-order language $\mathcal{L}$ which admits an infinite model $M$. Then for all cardinals $\kappa>|\mathcal{L}|$, there exists a first-order theory of size $\kappa$ which is elementary equivalent to $M$ (see, for instance,~\cite{leary_2015,marker_2002}).\\

One may heuristically understand this as saying that an infinite model of some theory can always be made more or less infinite, for instance, by extending it to a non-standard model. To understand the behaviour of such infinite models, which may, in general, be very badly behaved, it is valuable to understand categoricity:\\

\textbf{Definition:} \textit{An abstract elementary class $\mathfrak{K}$ is categorical in $\lambda$ if it has one and only one model of cardinality $\lambda$ up to isomorphism). The class of cardinals for which $\mathfrak{K}$ is categorical is called the categoricity spectrum of $\mathfrak{K}$ and is denoted $cat(\mathfrak{K})$.}\\

Categoricity allows us to compare the relative size of elements of some abstract elementary class, and indeed, there are several nice facts about categoricity which lend themselves nicely to the demarcation of distinct isomorphism classes. Notably, the following theorem is substantial~\cite{morley_1965}:\\

\textbf{Theorem (Morley):} \textit{Fix some countable first-order $\mathcal{L}$-theory $T$. Let $\lambda$ be a cardinal, and $I(\lambda,K)$ be the number of models in $K$ of cardinality $\lambda$ (up to isomorphism). Then if for some $\lambda>|T|+\aleph_0$, $I(\lambda,T)=1$, then $I(\mu,T)=1$ for all cardinals $\mu>|T|+\aleph_0$.}\\

That is, if the categoricity spectrum of some countable first-order theory includes a single model of some infinite cardinality larger than that of the language used, then indeed the categoricity spectrum of that theory extends to all such cardinals. Categoricity is the gateway to stable model theory which attempts to understand how well-constrained the infinite models of some theory are. If one wants to consider more specific properties unique to a models of a particular theory, the route to formalizing this is \textit{definability}, which tells you when the a particular theory admits the specification of particular subsets based on a precise set of properties which can be formulated in the language in question.\\

\textbf{Definitions:} \textit{Let $\mathfrak{A}$ be an $\mathcal{L}-structure$ with universe $A$ and let $B\subset A$. A set $X\subset A^n$ is $B$-definable if there exists some $\phi(v_1,...,v_n,w_1,...,w_n)$ such that $X=\{\overline{a}\in A^n|\phi(\overline{a},\overline{b}), \overline{b}\in B^m\}$.}\\

We can now say that classification theory deals with the problems of understanding isomorphism classes within abstract elementary classes of models, and frequently, the understanding of the isomorphism classes relies on what features are definable within a particular first-order theory and how the categoricity of the models behaves. The high-level goal of classification theory is to look at a collection of such classes $\{\mathfrak{K}\}_I$ and present clear mechanisms for demarcating the ``nice" or ``well-behaved" classes from those which are ``messy" or ``complicated". In this way, one might obtain a strong way to pin down the well-behaved classes, and prove substantive claims about all of these classes through more straight-forward case work. The most important result along these lines has been due to Shelah who determined the precise conditions for which an abstract elementary class admits a classification (or in more technical language, has a structure theory) of its isomorphism classes (the so-called dichotomy theorem)~\cite{shelah_2009,hodges_1987,shelah_1990}. This theorem states that:\\

\textbf{Theorem (Shelah):}\textit{An abstract elementary class $\mathfrak{K}$ has a structure theory if there is a class of subsets $\{A\}_I$ of the universe of each model in $K$ (called cardinal invariants), which are definable in each model, such that, for any two $\mathcal{M},\mathcal{N}\in K$, $\mathcal{M}\cong\mathcal{N}$ if and only if $|A_i^{\mathcal{N}}|=|A_i^{\mathcal{M}}|$ for all invariants $A_i$.}\\

The low-level goal of classification theory is, given a particular class $\mathfrak{K}$ which is known to admit a classification, to determine concretely all of the elements of $\mathfrak{K}$ up to isomorphism, or at least to determine a cardinality bound for how many such distinct elements there may be. The later part of this project was, again, solved by Shelah~\cite{shelah_2009,hodges_1987,shelah_1990}. In the operator algebra setting, once one has established that a particular kind of operator algebra with a first-order set of axioms in fact does admit a classification (a structure theory), then the important work is to determine every element of the class of models of those axioms up to isomorphism. The earliest results along these lines were the classification of UHF-algebras due to Glimm~\cite{rordam_2000,glimm_1960} and the classification of AF-algebras due to Elliott~\cite{rordam_2000,elliott_1976}.

Before proceeding to the operator algebra setting, to see a nice case of classification theory in action in a comfortable setting, consider vector spaces.\\

\textbf{Example:} \textit{Vector Spaces.} The class of all vector spaces $\mathcal{V}$ is the class of all models of the theory of vector spaces $T_{VS}$. One can show readily that $T_{VS}$ is a first-order theory~\cite{leary_2015}. Therefore, the class of vector spaces forms an abstract elementary class (ordered under inclusion, where all of the other more complicated requirements are satisfied more or less trivially). Given a vector space $V\in\mathcal{V}$, there is a definable subset of $V$ which is the set of all bases of $V$ (the relation under which this set is definable in the language of Marker~\cite{marker_2002}, is span). It is an elementary fact of linear algebra (and thus a theorem provable in the theory $T_{VS}$, and so true of all vector spaces) that the cardinality of every basis of a particular vector space is the same. Thus, the cardinality of every element of the set of bases of $V$ is constant. We call this cardinality the \textit{dimension} of $V$, and it is provable using elementary linear algebra that two finite-dimensional vector spaces are isomorphic if and only if their dimensions are the same. Therefore, the cardinal invariant for the structure theory of vector spaces is the dimension, which is the cardinality of any element of the definable subset of the universe of $V$ which is the set of bases for $V$.

\section{$C^*$-algebras and their K-Theory}

In the classification of $C^*$-algabras, K-theory is a valuable tool, as it allows one to construct a related structure to the algebra in question, namely the ordered $K_0$-group (as well as higher $K$-groups), which encodes lots of the information about the algebra while allowing questions to be posed in the context of abelian groups rather than operator algebras. Because of this, K-theory is helpful for the classification of $C^*$-algebras, as definable invariants for the algebra may be functorially pushed forward to the K-theory setting where they are more natural and easier to work with. In the special case of AF-algebras, the classification process is made very clear. After presenting the basic ideas pertinent to studying $C^*$-algebras and their K-theory, I walk through the examples of UHF-algebras and AF-algebras to understand their classification in model theoretic terms. I then discuss more general ideas presently being pursued for more general classifications. As such, I base the following on~\cite{rordam_2000}.\\

\textbf{Definition:} \textit{A $C^*$-algebra is a Banach algebra equipped with a conjugate-linear homomorphic isometric involution $x\mapsto x^*$ such that $||x^*x||=||x||||x^*||=||x||^2$.}\\

Broadly speaking, a $C^*$-algebra is an abstract representation of some algebra of bounded, self-adjoint, norm-closed operators on a Hilbert space. Indeed, the Gelfand-Naimark-Segal (GNS) construction tells us precisely that every abstract $C^*$-algebra admits such a Hilbert space representation. This notion of concrete representation is made stronger by the existence of the Gelfand transform, which acts as an isometric $^*$-isomorphism between elements of any given commutative, unital $C^*$-algebra and the set of bounded continuous functions on some particular compact Hausdorff space\footnote{In this sense, it may indeed be reasonable to think of the study of $C^*$-algebras as an abstract study of noncommutative topology.}. The study of specific operator algebras over particular Hilbert spaces together with their associated spectra is a significant area of study in its own right, but the abstract treatment of general $C^*$-algebras is likewise very rich, and particularly amenable for understanding the behaviour of such systems in the most unified, general setting.

Even though `zooming out' to the higher-level theory of abstract operator algebras allows one to forget about many of the precise details of the systems they are considering, $C^*$-algebras are, in their own right, highly complicated systems. As such, it is valuable to have tools to simplify problems pertaining to $C^*$-algebras further. One such method for gaining further understanding comes from K-theory.

When considering any operator algebra, the class of idempotent elements is of particular importance (one need only look towards von Neumann's formalism of quantum mechanics, or more modern POVM-based formalisms to see how practically powerful projections may be in the wild~\cite{nielsen_2010}). Specifically, projections in an algebra `pin down' much of the dimension structure of the algebra\footnote{For a particularly easy example, given the matrix algebras $M_n(\mathbb{C})$ and $M_k(\mathbb{C})$, one has that $M_n(\mathbb{C})\cong M_k(\mathbb{C})$ if and only if $n=Tr([I_n]_0)=Tr([I_k]_0)=k$ where $[I_n]_0$ and $[I_k]$ are the K-theoretic identities associated with each algebra (here, just the matrix identity up to Murray-von Neumann equivalence).} quite concretely in a precise sense using unitary equivalence and Murray-von Neumann equivalence:\\

\textbf{Definition:}\textit{Let $\mathcal{A}$ be a $C^*$-algebra. Let $p,q\in \mathcal{A}$. Then $p$ and $q$ are Murray-von Neumann equivalent (written $p\sim q$) if there exists some $v\in \mathcal{A}$ such that $p=v^*v$ and $q=vv^*$.}\\

It can readily be shown that if $p,q\in \mathcal{A}$ are projections which are homotopic or unitarily equivalent, then $p\sim q$. There is a canonical matrix algebra semigroup of projections associated with any $C^*$-algebra $\mathcal{A}$ given by:\\

\textbf{Definition:} \textit{Let $M_{n}(\mathcal{A})$ denote the algebra of $m\times m$ matrices with entries in the $C^*$-algebra $\mathcal{A}$. Let $\mathcal{P}_n(\mathcal{A})$ denote the class of projections in the matrix algebra $M_{n}(\mathcal{A})$. Then the projection semigroup is $\mathcal{P}_{\infty}(\mathcal{A}):=\cup_{n=1}^{\infty}\mathcal{P}_n(\mathcal{A})$.}\\

Addition in this semigroup is given by the usual direct sum $\oplus$ on matrices. There is an associated notion of Murray-von Neumann equivalence over $\mathcal{P}_{\infty}(\mathcal{A})$ given by:\\

\textbf{Definition:} \textit{Let $p\in \mathcal{P}_n(\mathcal{A})$ and $q\in\mathcal{P}_m(\mathcal{A})$. Then $p\sim_0 q$ if there is some $v\in M_{m,n}(\mathcal{A})$ such that $p=v^*v$ and $q=vv^*$.}\\

The triumph of K-theory is that, in light of this equivalence on the projection semigroup associated with a $C^*$-algebra $\mathcal{A}$, one can readily construct a \textit{functor} from the category of $C^*$-algebras to the category of ordered abelian groups, which preserves most of the relevant information about the structure of the algebra through its idemptotent elements. This functor is called the $K_0$ functor, and it has many beautiful features. The $K_0$ functor (which is a map which takes a $C^*$-algebra $\mathcal{A}$ to an ordered abelian group $K_0(\mathcal{A})$) is defined as follows:\\

\textbf{Definition:} \textit{Let $S$ be a semigroup. The Grothendieck group $G(S)$ is given by the quotient $G(S)=S\times S/\sim_G$ where $(x_1,y_1)\sim_G(x_2,y_2)$ if and only if there is some $z$ with $x_1+y_2+z=x_2+y_1+z$. The $K_0$ group associated with a $C^*$-algebra $\mathcal{A}$ is given by the Grothendieck group of the projection semigroup of $\mathcal{A}$. That is, $K_0(\mathcal{A})=G(P_\infty(\mathcal{A})/\sim_0)$.}\\

To better understand what $K_0(\mathcal{A})$ `looks like', the following results are quite valuable:\\

\textbf{Theorem:} \textit{ Let $[p]_0$ denote the $\sim_0$ equivalence class of some element $p\in\mathcal{P}_\infty(\mathcal{A})$. Then $K_0(\mathcal{A})=\{[p]_0-[q]_0:p,q\in\mathcal{P}_n(\mathcal{A}),n\in\mathbb{N}\}$. Moreover, $[p]_0+[q]_0=[p\oplus q]_0$. Finally, if $p$ and $q$ are homotopic, then $[p]_0=[q]_0$.}\\

Furthermore, the functoriality of $K_0$ tells us that it distributes in the natural sense over the composition of $^*$-homomorphisms and preserves the identity map. Additionally, if one has the categorial concept of inductive limits in place (described below), there is a natural sense in which $K_0$ is continuous.

With the notion of $K_0$ in place, we may begin to talk about specific classes of $C^*$-algebras in terms of their K-theory. The most elementary class of $C^*$-algebras is that of \textit{finite} $C^*$-algebras (each of which is isomorphic to the matrix algebra $M_n(\mathbb{C})$ for some $n$). These algebras are not terribly exciting, but one can form many highly non-trivial classes of $C^*$-algebras out of them using inductive limits. In their general construction, inductive limits act like limits of sequences of objects in some category which are linked to each other by connecting arrows; they generalize limits to more abstract categories. Specifically:\\

\textbf{Definition:} \textit{Let $C$ be some category. Let $\{O_i,\phi_{ij}\}$ be a sequence of objects $O_i$ in $C$ and arrows $\phi_{ij}:O_i\to O_{j}$ ($i<j$). Then the inductive limit $\lim_{\to}\{O_i,\phi_{ij}\}$ is the object $O_\infty$, which has a canonical set of maps $\phi_i:O_i\to O_\infty$ such that $O_\infty=\cup_{i=1}^\infty \phi_i(O_i)$ and where $\phi_j=\phi_i\circ\phi_{ij}$ (i.e. the below diagram commutes).}\\

\begin{center}
\begin{tikzcd}
O_1 \arrow[rrr, bend left, "\phi_1"] \arrow[r, "\phi_{12}"] & O_2 \arrow[rr, bend right, "\phi_2"]\arrow[r,"\phi_{23}"]&\arrow[r]...&O_\infty
\end{tikzcd}
\end{center}

It is important to note that not all categories admit inductive limits (for instance, the category of sets). When we take the category in question to be the category of $C^*$-algebras, we can form inductive limits (taking the connecting maps to be $^*$-homomorphisms), and it turns out that even the inductive limits of finite $C^*$-algebras can be quite fascinating objects to study. For our interests here, this setting also provides us with a straightforward and historically significant instance of classification theory found in the wild. Specifically, the study of UHF-algebras and AF-algebras (which are inductive limits of finite $C^*$-algebras) provide an interesting setting to connect Shelah's abstract classification theory to concrete and relevant classification programs.

\section{The Classification of $C^*$-algebras}

To proceed, we begin with the following definition of approximately finite (AF) algebras:\\

\textbf{Definition:} \textit{An AF-algebra is a $C^*$-algebra which is isomorphic to the inductive limit of a sequence of finite $C^*$-algebras.}\\

The first example of classification theory as it pertains to $C^*$-algebras which makes clear how these ideas are all inter-connected in a concrete way is the classification of uniformly hyper-finite (UHF) algebras due to Glimm~\cite{glimm_1960}. An algebra $\mathcal{A}$ is a UHF-algebra if it is an AF-algebra where the connecting $^*$-homomorphisms are unit preserving. UHF-algebras may be classified using supernatural numbers which are defined as~\cite{rordam_2000}:\\

\textbf{Definition:} \textit{Let $\{n_i\}_{i\in\mathbb{N}}$ be a (finite or countable) set of cardinals where $n_i\leq \aleph_0$, and let $n$ formally be expressed as the product $n=\Pi_{i\in\mathbb{N}}p^{n_i}_i$ where $p_i$ is the $i$th prime number. Any such $n$ is a supernatural number.}\\

It will be shown that $\{n_i\}_{i\in\mathbb{N}}$ is the set of cardinal invariants of models of $T_{UHF}$\footnote{It should be noted that, in general, $C^*$-algebras are not universally first-order axiomatizable, and even when they are, they generally require an uncountable language, making certain results, such as L\"owenheim-Skolem results, more subtle. However, it can be shown that abelian algebras, and unital finite $C^*$-algebras (as well as many other elementary kinds of algebras) are first-order axiomatizable, and so we do not get into trouble here~\cite{farah_2018}.}, the first-order theory of UHF algebras. Without getting our hands too dirty, stated without proof, we have the following theorems~\cite{rordam_2000}:\\

\textbf{Theorem:} \textit{Let $G\subset(\mathbb{Q},+)$ be a subgroup containing 1. Then there is a supernatural number $n$ associated with $G$ such that $G=Q(n)$, the subgroup composed of elements of $\mathbb{Q}$ where the denominator is a factor of $n$.}\\

\textbf{Theorem:} \textit{Let $\mathcal{A}$ be a UHF algebra. Then $(K_0(\mathcal{A}),[1_\mathcal{A}]_0)\cong (Q(n),1)$ for some supernatural number $n$.}\\

From this, we see that every UHF-algebra has an associated supernatural number $n$, and so to every UHF-algebra, a set of cardinal invariants $\{n_i\}_{i\in\mathbb{N}}$ may be assigned. This is a step towards classification, however, we do not yet know that this associated supernatural number determines a given UHF-algebra \textit{uniquely} up to isomorphism; \textit{a priori}, there may exist two non-isomorphic UHF-algebras (models of $T_{UHF}$) with the same associated supernatural number. However, the major result of Glimm says that the associated supernatural number of a UHF-algebra is in fact enough to determine that UHF-algebra up to isomorphism. Specifically, we have:\\

\textbf{Theorem: (Glimm):} \textit{If $\mathcal{A}$ and $\mathcal{A}'$ are two UHF algebras with associated supernatural numbers $n$ and $n'$ respectively, then $\mathcal{A}\cong\mathcal{A}'$ if and only if $n=n'$.}\\

In the model-theoretic setting, if $n=n'$, since prime factorizations are unique, we see that $n_i=n_i'$ for all elements of $\{n_i\}_{i\in\mathbb{N}}$ and $\{n'_i\}_{i\in\mathbb{N}}$, and so we see that $\mathcal{A}\cong\mathcal{A}'$ if and only if $|n_i|=|n'_i|$ for all $i\in\mathbb{N}$ and so $\{n_i\}_{i\in\mathbb{N}}$ is a (countable) set of cardinal invariants which constitute the structure theory for $T_{UHF}$. Therefore, we see a nice reflection of Shelah's classification theory in Glimm's theory for UHF-algebras.

In light of the fact that $(K_0(\mathcal{A}),[1_\mathcal{A}]_0)\cong (Q(n),1)$ from above, we see that in fact we may restate this classification theorem to be in the language of K-theory. Specifically, $\mathcal{A}\cong\mathcal{A}'$ if and only if $(K_0(\mathcal{A}),[1_\mathcal{A}]_0)\cong(K_0(\mathcal{A'}),[1_\mathcal{A'}]_0)$, because the $K$ groups are determined by the associated supernatural numbers.

Further connecting the classification of UHF-algebras to classification theory in general, if we recall the low-level goal of determining the cardinality of the set of isomorphism classes of models for a particular theory ($T_{UHF}$ in this case), we see that there is a converse theorem which comes in handy. It is shown in~\cite{rordam_2000} that for every supernatural number $n$, there is a UHF-algebra with associated supernatural number $n$. A supernatural number is exactly specified by the set $\{n_i\}_{i\in\mathbb{N}}$ (noting that we have here ordered the $n_i$'s using the same usual ordering of the corresponding prime numbers, and that prime factorization is unique). Since there are $\aleph_0$ elements of $\{n_i\}_{i\in\mathbb{N}}$, and since each $n_i$ may take on $\aleph_0$ possible values, there are $\aleph_0^{\aleph_0}$ many supernatural numbers. One may show that $\aleph_0^{\aleph_0}=2^{\aleph_0}=\aleph_1=|\mathbb{R}|$, and so there is a (small) uncountable set of distinct isomorphism classes of UHF algebras. I now look at unital AF-algebras in general.

Let $T_{AF}$ be the first-order theory of unital AF-algebras (everything to follow may be generalized to the non-unital case~\cite{rordam_2000,elliott_1976}, but the details involved do not lead to a clearer understanding of the general classification method). The question which must be answered to classify AF-algebras in the sense of Shelah is: What set of cardinals $\{\lambda_i\}_{i\in I}$ is definable in any model $\mathcal{A}\vDash T_{AF}$ such that $\{\lambda_i\}_{i\in I}=\{\lambda'_i\}_{i\in I}$ if and only if $\mathcal{A}\cong\mathcal{A}'$? The theorem of Elliott~\cite{elliott_1976} of course indicates how this may be answered:\\

\textbf{Theorem (Elliott):} \textit{Two unital AF-algebras $\mathcal{A}$ and $\mathcal{B}$ are isomorphic if and only if $(K_0(\mathcal{A}),K_0(\mathcal{A})^+,[1_\mathcal{A}]_0)\cong(K_0(\mathcal{B}),K_0(\mathcal{B})^+,[1_\mathcal{B}]_0)$.}\\

Thus, the ordered $K_0$ group with its unit is an algebraic invariant which is sufficient to classify AF-algebras. But how do we interpret this in terms of sets of cardinals? The following result, found in~\cite{rordam_2000}, brings clarity to the matter:\\

\textbf{Theorem:} \textit{Let $\mathcal{A}$ be an AF algebra. Then $(K_0(\mathcal{A}),K_0(\mathcal{A})^+)$ is a dimension group. Conversely, for every dimension group $(G,G^+)$, there exists an AF-algebra $\mathcal{A}$ such that $(G,G^+)\cong(K_0(\mathcal{A}),K_0(\mathcal{A})^+)$.}\\

The notion of dimension is valuable, because it allows us to tease out the cardinals from the discussion. As is explained clearly by Effros et al.~\cite{effros_1980}, `dimension' may take on many different meanings. In the context of matrix algebras, dimension refers to the rank of the projections in the algebra. In the von Neumman algebra setting, dimension becomes a non-negative real value. In the $C^*$-algebra setting however, dimension refers to a value assigned to Murray-von Neumman equivalence classes of projections in the $K_0$-group of the algebra; dimension then refers to values in some ordered abelian group. In this way, we see why the term dimension group is fitting in the AF-algebra setting; the dimension group is precisely the $K_0$ group in which elements of the algebra may have their `dimension' defined.\\

\begin{figure}[ht]
    \centering
    \begin{tikzcd}
    \cdot^1 \arrow[d, equal]\\
    \cdot^2 \arrow[d, equal]\\
    \cdot^4 \arrow[d, equal]\\
    \cdot^8 \arrow[d, equal]\\
    \vdots
    \end{tikzcd}
    \caption{The Bratteli diagram for the inductive sequence of finite-dimensional $C^*$-algebras given by $\mathcal{A}_n\cong M_{2^n}(\mathbb{C})$, with the connecting maps given by $\phi_{ii+1}(x)=diag(x,x)$. We see that the inductive limit of this sequence is isomorphic to $M_{2^\infty}(\mathbb{C})$. It may be shown that this is also the (unique) $C^*$-completion of the canonical anti-commutation relation (CAR) algebra over a separable infinite-dimensional Hilbert space, an algebra which is very important in mathematical physics, the study of quantum field theories and quantum statistical mechanics in particular (anything which makes use of antisymmetric Fock space).}
    \label{fig:CAR}
\end{figure}
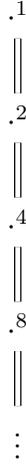

If we wish to understand how Elliott's classification theorem for AF-algebras may be expressed in terms of cardinal invariants for the abstract elementary class of models of the first-order theory for AF-algebras, it is valuable to follow Bratteli~\cite{bratteli_1972} in his diagrammatic representation of AF-algebras using infinite directed graphs, and read off the cardinal invariants from this graph in terms of counting the vertices and edges in a particular way. Bratteli diagram for an AF-algebra $\mathcal{A}\cong \lim_{\to}\{\mathcal{A}_i,\phi_{ij}\}$ are constructed thusly: For each $\mathcal{A}_i$, there will be a row of the Bratteli diagram consisting of a collection of points. Each $\mathcal{A}_i$ is isomorphic to some finite factor space of finite-dimensional complex matrix algebras. That is, $\mathcal{A}_i\cong \bigoplus_{k=1}^{s}M_{k_i}(\mathbb{C})$. Thus, in the $i$th row of the Bratteli diagram, draw one point for each factor $M_{k_i}(\mathbb{C})$ of $\mathcal{A}_i$, and label it by the integer $k_i$. Given the $i$th row of a Bratteli diagram, the $i+1$th row is obtained by likewise drawing one point for each factor space of $\mathcal{A}_{i+1}$ (in its finite-dimensional complex matrix algebra representation, for instance). Then the $i$th row is connected to the $i+1$th row by drawing edges connecting the points in the $i$th row to points in the $i+1$th row precisely when a given point in the $i$th row is mapped to the corresponding point in the $i+1$th row under the connecting map $\phi_{ii+1}$. The number of lines connecting two points denotes the \textit{multiplicity} of the map connecting the two points; if $\phi_{ij}:\mathcal{A}_i\to\mathcal{A}_j$ is a $^*$-homomorphism, then the multiplicity of $\phi_{ij}$ is defined to be $Tr(\phi_{ij}(p))/Tr(p)$ for any non-zero projection $p$ (it is easy to show that this quantity is independent of the choice of $p$). It turns out this graph representation completely characterizes the inductive limit and thus completely characterizes up to isomorphism the underlying AF-algebra. Two examples of Bratteli diagrams (taken from~\cite{rordam_2000}) may be found in figures \ref{fig:Fib} and \ref{fig:CAR}.\\

\begin{figure}[h]
    \centering
    \begin{tikzcd}
    \cdot^1 \arrow[d, dash] \arrow[rd, dash] & \cdot^1 \arrow[ld, dash]\\
    \cdot^2 \arrow[d, dash] \arrow[rd, dash] & \cdot^1 \arrow[ld, dash]\\
    \cdot^3 \arrow[d, dash] \arrow[rd, dash] & \cdot^2 \arrow[ld, dash]\\
    \vdots & \vdots
    \end{tikzcd}
    \caption{The Bratteli diagram for the inductive sequence of finite-dimensional $C^*$-algebras constructed by taking $\mathcal{A}_n=M_{f_n}(\mathbb{C})\oplus M_{f_{n-1}(\mathbb{C})}$ where $f_0,f_1,f_2,\dots$ are the usual Fibonacci numbers, with connecting maps given by: $(x,y)\mapsto\big(\usebox{\smlmat},x\big)$.}
    \label{fig:Fib}
\end{figure}
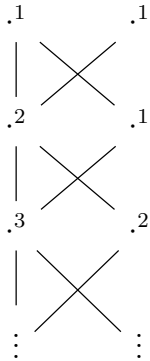

The ordered abelian $K_0$ group (or the dimension group, as one chooses) associated with an AF-algebra, as Elliott's theorem showed, comes from the sequence of the $K_0$-groups of the terms in the inductive sequence of finite $C^*$-algebras which induce the AF-algebra in question. As such, to determine the algebraically invariant dimension group of an AF-algebra, one may look at the dimension groups associated with each algebra in the sequence. But this is effectively the content of a Bratteli diagram! Bratteli diagrams show us, at each level, the dimension of each finite matrix sub-algebra (using the unusual notion of dimension, which may readily be cast as a cardinal) together with the way in which these finite matrix algebras are `stitched' together via their connecting maps, as well as the multiplicity of these connecting maps (which therefore determines how the rank of projections changes when pushed forward to the next row of the Bratteli diagram). The `bottom row' of a Bratteli diagram thus determines exactly the dimension group of the inductive limit of the sequence, which is therefore isomorphic to the $K_0$ group of the AF-algebra in question.

The $i$th row of a Bratteli diagram may be completely specified by the number associated to each point, together with the edges connecting it to the previous row. Thus, the $i$th row may be specified completely by the set of cardinals $\{n_i,\{d_{i,j}\}_{j=1}^{n_i},\{m_{(i-1,s)\to(i,r)}\}_{s,r=1}^{s=n_{i-1},r=n_i}\}$ where $n_i$ is the number of elements of the $i$th row,  $d_{i,j}$ is the dimension associated with the point in the $j$th position (from left to right) of the $i$th row, and $\{m_{(i-1,s)\to(i,r)}\}_{s,r=1}^{s=n_{i-1},r=n_i}$ denotes the multiplicities of each edge connecting the $s$th point in the $i-1$th row to the $r$th point in the $i$th row (where each $m_{(0,s)\to(1,r)}=0$ for all $s$ and $r$ since there are no connecting maps ending in the first row). These values are all finite integers for any given row. The `bottom line' (i.e. the $i$th row in the limit when $i\to\aleph_0$, the countable infinity) may likewise be specified by this set of cardinals.\footnote{For the set theorist, this is an terribly dull set of cardinals, as all of them are finite or countably infinite. However, boring or not, cardinals they are, and so they are satisfactory for us here.} The Bratteli diagram encodes all of the information about the $K_0$ group of the AF-algebra in question, and so it is an equivalent description of the invariant necessary to classify AF-algebras, although it is more conducive to being read as a sequence of cardinals, formally putting us in a position to claim, in Shelah's language, that models of the theory $T_{AF}$ are classified up to isomorphism by a set of cardinal invariants.

UHF and AF-algebras are not the only classes of $C^*$-algebras to be classified, nor are they the only classes of $C^*$-algebras defined in terms of inductive limits to be classified. However, noting that, in the UHF-algebra case, the invariant required could be equivalently expressed as the double $(K_0(\mathcal{A}),[1_\mathcal{A}]_0)$, and in the more general unital AF-algebra case, the invariant needed was only the slightly more sophisticated triple $(K_0(\mathcal{A}),K_0(\mathcal{A})^+,[1_\mathcal{A}]_0)$, one might be inclined to state a more general thesis, namely that all $C^*$-algebras of some form may be classified by a more general K-theoretic invariant. This is the present program of classification, and one of the conjectured invariants is the so-called Elliott invariant~\cite{farah_2019}:

$$Ell(\mathcal{A})=(K_0(\mathcal{A}),K_0(\mathcal{A})^+,[1_\mathcal{A}]_0,K_1(\mathcal{A}),T(\mathcal{A}),r_\mathcal{A})$$

Here, $K_1(\mathcal{A})$ is the standard $K_1$ group (which is isomorphic to the $K_0$ group of $S\mathcal{A}$, the topological suspension of $\mathcal{A}$), $T(\mathcal{A})$ is the tracial space of $\mathcal{A}$ (i.e. the set of all tracial functionals in the dual space of $\mathcal{A}$), and $R_\mathcal{A}$ is the so-called coupling map which associates states with traces. The goal of the classification program was to classify all separable, unital, nuclear simple $C^*$-algebras via this invariant. That is, ideally, one would end up with a theorem which states that, if $\mathcal{A}$ and $\mathcal{B}$ are two such algebras, then $\mathcal{A}\cong\mathcal{B}$ if and only if $Ell(\mathcal{A})\cong Ell(\mathcal{B})$. Unfortunately, there have been counter-examples preventing such a complete functorial classification, but it has held for many important cases, such as the aforementioned algebras, as well as irrational rotation algebras, and AT-algebras in general. There are presently new candidate invariants being proposed. However, these classification programs have become extremely complicated and sophisticated, and a complete functorial classification may be far away; one of the major problems in this field currently is to try and understand precisely how complex the isomorphism maps between $C^*$-algebras may in principle be, thereby giving insight into the complexity of their classification~\cite{farah_2019}.

\section{Conclusions}

I here presented a general overview of Shelah's classification theory, and the theory of cardinal invariants. After a cursory glance at the basics of K-theory for $C^*$-algebras, I then recalled several known results, due primarily to Glimm and Elliott, about the classification of various classes of $C^*$-algebras using particular K-theoretic invariants. I demonstrated how these invariants may be understood as cardinal invariants in terms of sequences defining supernatural numbers and sequences defining states of Bratteli diagrams. With this understanding in place, I described the goals of future classification. The topic of $C^*$-algebra classification has, in recent years, become deeply rooted in model theory and mathematical logic, and so framing the discussion in terms of Shelah's classification theory reorients the setting of discussion to be more amenable to such discourse. In this way, we see that the marriage of model theory and operator algebra theory is, in a sense, preordained by the sophistication of the field, and will likely serve as an important tool for future investigation.

\bibliographystyle{unsrt}
\bibliography{classification_theory}

\end{document}